\documentclass[12pt]{amsart}
\usepackage{mathpazo,fullpage,relsize}
\usepackage{latexsym,amscd,amssymb}
\newtheorem{thm}{Theorem}
\newtheorem{prop}[thm]{Proposition}
\newtheorem{lem}[thm]{Lemma}

\theoremstyle{remark}
\newtheorem{rem}[thm]{Remark}

\theoremstyle{definition}
\newtheorem{defn}[thm]{Definition}

\newcommand{\C}{\mathbb{ C}}

\newcommand{\lra}{\longrightarrow}

\newcommand{\Q}{\mathbb{ Q}}

\newcommand{\HH}{\mathbb{ H}}

\newcommand{\Z}{\mathbb{ Z}}

\makeatletter
\@namedef{subjclassname@2020}{%
  \textup{2020} Mathematics Subject Classification}
\makeatother

\title[Unbounded Pontryagin numbers on nonnegatively curved spin manifolds]{Unbounded Pontryagin numbers on\\
nonnegatively curved spin manifolds}
\author{E.~Hsiao}
\address{Mathematisches Institut, {\smaller LMU} M\"unchen,
Theresienstr.~39, 80333~M\"unchen, Germany}
\email{enya0034@gmail.com}
\author{D.~Kotschick}
\address{Mathematisches Institut, {\smaller LMU} M\"unchen,
Theresienstr.~39, 80333~M\"unchen, Germany}
\email{dieter@math.lmu.de}

\date{November 11, 2021, revised April 9, 2022; \copyright{\ E.~Hsiao and D.~Kotschick 2021}}
\subjclass[2020]{primary 53C20, 57R20; secondary 53C23, 53C27, 57R75}

\begin{document}

\begin{abstract}
We prove that any rational linear combination of Pontryagin numbers that does not factor through the 
universal elliptic genus is unbounded on connected closed spin manifolds of nonnegative sectional curvature. 
\end{abstract}

\maketitle


\section{Introduction}

One of the most important results in the study of Riemannian manifolds with nonnegative curvature 
is Gromov's Betti number theorem~\cite{G}, which gives, in every dimension, a universal upper bound 
for the Betti numbers of connected nonnegatively curved manifolds, thus bounding the size of such
manifolds. As a corollary, any invariant that can be bounded in terms of Betti numbers is bounded 
on nonnegatively curved manifolds. This applies in particular to the signature or $L$-genus. In sharp contrast
with this, all other linear combinations of Pontryagin numbers are not bounded on nonnegatively 
curved manifolds. It was proved in~\cite{K} that, up to taking multiples, the $L$-genus is characterised 
among all linear combinations of Pontryagin numbers by its boundedness on  connected manifolds of 
nonnegative curvature. 

It is clear that the result of ~\cite{K} cannot remain true if one restricts the discussion to spin manifolds, 
since by the Lichnerowicz argument~\cite{L} the $\hat{A}$-genus vanishes on nonnegatively curved spin 
manifolds, and, starting in dimension $8$, it is linearly independent of the $L$-genus. However, both the 
$L$-genus and the $\hat{A}$-genus are specialisations of Ochanine's elliptic genus~\cite{O}, and this led Herrmann 
and Weisskopf~\cite{HW} to ask whether linear combinations of Pontryagin numbers which do not arise 
from the elliptic genus are unbounded on connected spin manifolds of nonnegative curvature. 
It was shown in~\cite{HW} that the answer to this question is positive up to dimension $20$, using 
various case by case calculations on suitably chosen examples, which however do not provide a
path towards extending the result to arbitrary dimensions. 
The question posed by Herrmann and Weiskopf is part of a general philosophy trying to 
find index-theoretic obstructions for the existence of positively or nonnegatively curved Riemannian
metrics, and suggests that the only such obstructions should come from the elliptic genus, cf. Dessai~\cite{D}.

In this paper we prove a result in the spirit of~\cite[Theorem~1]{K} for spin manifolds, 
which in particular provides a positive answer to Question~1.2 in~\cite{HW}:
\begin{thm}\label{t:main}
Any rational linear combination of Pontryagin numbers that does not factor through the elliptic genus 
is unbounded on connected closed oriented spin manifolds of nonnegative sectional curvature. 
\end{thm}
The proof uses the structure of the rational spin bordism ring. We construct sequences of ring 
generators with the property that each generator of dimension $\geq 12$ belongs to a family of 
possible choices on which a certain indecomposabe Pontryagin number is unbounded. Moreover, 
with the same dimension assumption, these generators admit metrics of nonnegative sectional 
curvature and are in the kernel of the elliptic genus. In terms of dimensions, this is the best one 
can hope for, since in dimensions up to $8$ all Pontryagin numbers factor through the elliptic genus. 
In dimension $4$, because of the Lichnerowicz theorem~\cite{L}, it is not possible to choose a spin 
generator with nonnegative curvature for the bordism group. This makes the proof of  Theorem~\ref{t:main}
much more complicated than the proof in~\cite{K}, since one has to work around the lack of a suitable 
$4$-dimensional generator. This is achieved by uncovering certain polynomial relations between different
families of generators for the bordism ring, see Proposition~\ref{noK3} in Section~\ref{s:bordism} below.

\section{A family of projective bundles over complex projective spaces}

For integers $c$ denote by $H^c$  the tensor powers of the hyperplane bundle
$H\lra \C P^n$. (For negative $c$ these are the $\vert c\vert$-fold tensor products of 
the dual line bundle.) Consider the complex vector bundle $E_{c} \longrightarrow \C P^n$ 
of rank $k+1$ of the form $E_c=H^c\oplus \underline{\C }^{k}$, 
and let $X_n^k(c)=P(E_c)$ denote its projectivization. This is a complex manifold of 
real dimension $2(n+k)$. To have non-trivial Pontryagin numbers we need manifolds of 
real dimension divisible by $4$, and so we will take $n$ and $k$ of the same parity.

The total Chern class of $E_c$ is $c(E_c)=1+c\cdot x$, where $x\in H^2(\C P^n;\Z)$ is the positive generator.
Therefore, the Leray--Hirsch theorem immediately gives the following statement.
\begin{lem}\label{cohomology}
The cohomology ring of $X_n^k(c)$ is generated by two classes $x,y\in H^2(X_n^k(c))$ subject to the 
relations 
\begin{equation}
x^{n+1}=0 \ , \quad \quad y^{k+1}+cxy^{k}=0 \ .
\end{equation}
\end{lem}
Here $x$ is the generator coming from the base $ \C P^n$.
Next we compute the characteristic classes of the $X_n^k(c)$.
\begin{lem}
The total Chern class of $X_n^k(c)$ is given by
\begin{align}
\label{total chern class}
    c(X_n^k(c))=(1+x)^{n+1} (1+y)^{k}(1+y+cx) \ .
\end{align}
Therefore, its total Pontryagin class is 
\begin{align}
\label{total pontryagin class}
    p(X_n^k(c))=(1+x^2)^{n+1} (1+y^2)^{k}(1+(y+cx)^2) \ .
\end{align}
\end{lem}
\begin{proof} 
We have the decomposition $TX_n^k(c)= T\pi \oplus \pi^*T\C P^n$, where $T\pi$ is the tangent bundle 
along the fibers, and so
$$
c(X_n^k(c)) = c(T\pi)\cdot \pi^* c( T\C P^n) 
$$
with the latter factor being $\pi^*c(\C P^n)=(1+x)^{n+1}$. 
It remains to compute $c(T\pi)$. 
We have the relative Euler sequence 
$$
0 \to L^{-1} \to \pi^*E_c \to L^{-1}\otimes T\pi \to 0 \ , 
$$
where $L$ is the fiberwise hyperplane bundle with $c(L)=1+y$. Tensoring by $L$ yields an isomorphism of bundles 
$L\otimes \pi^*E_c \cong T\pi\oplus \underline{\C }$, and this shows 
$$
c (T\pi)=c (L\otimes \pi^*E_c ) = (1+y)^{k}(1+y+cx) \ .
$$
This proves the formula for the total Chern class, which in turn gives the one for the total Pontryagin class.
\end{proof}

The following lemma tells us under which conditions $X_n^k(c)$ is spin.
\begin{lem}\label{spinlemma}
Assume that $n$ and $k$ have the same parity. Then $X_n^k(c)$ is spin
if and only if $k$ and $n$ are odd and $c$ is even. 
\end{lem}
\begin{proof} 
Recall that a complex manifold is spin if and only if its first Chern class is divisible by $2$ in integral
cohomology. For $X_n^k(c)$ the formula~\eqref{total chern class} gives us
$$
c_1(X_n^k(c))=(n+1)x +ky +(y+cx)=(k+1)y+(n+1+c)x \ .
$$
This is divisible by $2$ if and only if both $k+1$ and $n+1+c$ are even. Since we assumed that $k$ and $n$ have 
the same parity this happens if and only if $k$ and $n$ are odd and $c$ is even.
\end{proof}

In later sections we will always assume that the conditions in this Lemma are satisfied, so that we are dealing
with spin manifolds of real dimension $4m$, where $2m=k+n$. To end this section, we note that the manifolds
$X_n^k(c)$ are nonnegatively curved.

\begin{lem}
Every $X_n^k(c)$ admits a Riemannian metric of nonnegative sectional curvature.
\end{lem}
\begin{proof}
Since $X_n^k(0)$ is a product of complex projective spaces, we may assume $c\neq 0$.
The vector bundle $E_c=H^c \oplus \underline{\C }^k$ over $\C P^n$ has structure group $U(1)=S^1$,
and is therefore associated to an $S^1$-bundle $P_c\lra\C P^n$ with Euler class $c\cdot x$. 
The total space of this circle bundle is the lens space
$S^{2n+1}/\Z_{\vert c\vert }$, carrying an $S^1$-invariant metric of constant positive curvature.
Now $X_n^k(c)$ has the form $(P_c\times \C P^k)/S^1$, where $S^1$ acts freely by isometries
of the product metric formed by the positively curved metric on the lens space and the 
Fubini--Study metric on $\C P^{k}$, since $S^1$ acts on $\C P^{k}$ via the appropriate 
inclusion $S^1\hookrightarrow U(k+1)$.
The non-decreasing property of curvature in submersions implies that the induced metric on 
$X_n^k(c)$ is nonnegatively curved.
\end{proof}

\section{Computations of Pontryagin numbers}

We now prove some results about the  Pontryagin numbers of $X_n^k(c)$ in the 
case where $k$ and $n$ are both odd. Here by Pontryagin number we mean not
just the evaluation of monomials in the Pontryagin classes, but any rational linear 
combination of such evaluations.
The general shape of these numbers is as follows.

\begin{prop}\label{prop:odd}
As a function of $c$, every Pontryagin number of $X_n^k(c)$ is an odd polynomial 
of degree at most $n$. 
\end{prop}
Recall that a polynomial is odd if it involves only odd powers of the variable $c$.
We will see in Proposition~\ref{Thom number lemma} below that the maximal degree $n$ does occur.
\begin{proof}
Recall from~\eqref{total pontryagin class} that the total Pontryagin class is given by 
$$
p(X_n^k(c))=(1+y^2)^{k}(1+x^2)^{n+1}(1+(y+cx)^2) \ .
$$
Therefore, every Pontryagin number is given by the evaluation of a homogeneous
polynomial  of degree $m=\frac{1}{2}(k+n)$ in the variables $y^2$, $x^2$ and $(y+cx)^2$.
In other words, we are looking into evaluating linear combinations of monomials
of the form
$$
y^{2a}\cdot x^{2b}\cdot (y+cx)^{2(m-a-b)} 
$$
with nonnegative exponents. Expanding the third factor with the binomial theorem
we find a linear combination of the monomials
$$
y^{2a+i}\cdot x^{2m-2a-i}\cdot c^{2(m-a-b)-i} \ \ \ \textrm{with} \ \ \ 0\leq i\leq 2(m-a-b) \ .
$$
Note that $x$ and $y$ are cohomology classes of degree $2$, and $c$ is an integer.
Setting $j=2m-2a-i$, these monomials become
$$
y^{2m-j}\cdot x^{j}\cdot c^{j-2b} \ \ \ \textrm{with} \ \ \ 0\leq j\leq n \ .
$$

Using the cohomology relation from Lemma~\ref{cohomology} repeatedly, we find
$$
y^{2m-j}\cdot x^{j}\cdot c^{j-2b} = (-1)^{n-j} y^k\cdot x^n \cdot c^{n-2b} \ .
$$
Here $y^k\cdot x^n$ is the generator of the top-degree cohomology of $X_n^k(c)$.
Since $n$ is odd, the exponent $n-2b$ of $c$ is always odd, and this finally 
shows that any Pontryagin number of $X_n^k(c)$ is a rational linear combination
of terms which contain only odd powers of $c$. Moreover, the exponent $n-2b$ of $c$ is 
bounded above by $n$ since $b$ is nonnegative.
\end{proof}

For any closed oriented manifold $M$ of dimension $4m$ with total Pontryagin class
$$
p(TM)=\mathlarger{\prod}_i (1+y_i^2) 
$$
the Milnor--Thom number $s_n(M)$  is defined by 
$$
s_m(M)  = \mathlarger{\sum}_i \bigl< y_i^{2m},[M]\bigr> \ . 
$$
The splitting principle implies that this is a Pontryagin number.
The significance of $s_m$ is that its non-vanishing on $M$ is equivalent
to $M$ being a generator of the rational bordism ring, see Section~\ref{s:bordism} below.

\begin{prop} 
\label{Thom number lemma}
Assume that $k$ and $n$ are odd, and let $m=\frac{1}{2}(k+n)$.
Then the Milnor--Thom number of $X_n^k(c)$ is given by
\begin{align}
\label{Thom number formula}
    s_{m}(X_n^k(c))=c^n\left[\binom{k+n-1}{n}-k\right] \ . 
\end{align}
In particular $s_{m}(X_n^k(c))\neq 0$ whenever $k\geq 3$, $n\geq 3$ and $c\neq 0$.
\end{prop}
\begin{proof} 
This formula is a special case of a calculation of Schreieder~\cite{Sch}, who considered 
arbitrary projectivizations. Nevertheless, we include a direct proof, which is no more complicated than
the explanation of how to extract what we need from~\cite[Lemma~2.4.]{Sch}.

Using the formula~\eqref{total pontryagin class}, the definition of the Milnor--Thom number 
gives
\begin{equation*}
    \begin{aligned}
    s_{m}(X_n^k(c)) &=\bigl< k y^{n+k}
    +(n+1)x^{n+k}+ (y+cx)^{n+k} , [X_n^k(c)]\bigr> \ . 
    \end{aligned}
\end{equation*}
The term $(n+1)x^{n+k}$ vanishes since $x^{n+1}=0$. We expand $(y+cx)^{n+k} $ using 
the binomial formula, and drop all terms where the exponent of $x$ is $>n$.
Finally we trade all terms with an exponent of $y$ that is larger than
$k$ using the relation $y^{k+1}=-cxy^k$ from Lemma~\ref{cohomology} repeatedly. 
This leads to
\begin{align}
    k y^{n+k}+ (y+cx)^{n+k}&= k y^{n+k}+ \sum_{i=0}^{n} \binom{n+k}{i}c^i x^i y^{n+k-i}  \nonumber \\
    &= (-1)^n k c^n x^n y^{k}+\sum_{i=0}^{n}(-1)^{n-i} \binom{n+k}{i}c^n x^n y^{k}  \nonumber \\
    &=(-1)^n c^n\left[k+\sum_{i=0}^{n}(-1)^{i} \binom{n+k}{i}\right] x^n y^k\ . \nonumber
\end{align}
Since $n$ is odd, the term $(-1)^n=-1$. Further,  $x^n y^k$ evaluates as $1$ on the fundamental 
class, and therefore
\begin{equation*}
    \begin{aligned}
    s_{m}(X_n^k(c)) &= - c^n\left[k+\sum_{i=0}^{n}(-1)^{i} \binom{n+k}{i}\right] \ . 
    \end{aligned}
\end{equation*}
Now replace each binomial coefficient with $i > 0$ in the sum using the recursion
\begin{align}
\label{recursive identity}
    \binom{n+k}{i}=\binom{n+k-1}{i-1}+\binom{n+k-1}{i}.
\end{align}
Since in the resulting sum almost all binomial coefficients appear twice with opposite
signs, they cancel in pairs, except for the last summand
$\binom{n+k-1}{n}$. This finally gives
  \begin{equation*}
    \begin{aligned}
    s_{m}(X_n^k(c)) &= - c^n\left[k-\binom{n+k-1}{n}\right] \ 
    \end{aligned}
\end{equation*}
as claimed.
\end{proof}

We will need to consider another special Pontryagin number, defined as follows.
\begin{defn}
For a closed oriented manifold $M$ of dimension $4m$ the $q$-number is defined by 
$$
q_m(M)=\mathlarger{\sum}_{i} \mathlarger{\sum}_{j\neq i} \bigl<  y_i^2  y_j^{2m-2} , [M]\bigr> \ ,
$$
where the $y_i^2$ are the Pontryagin roots of $M$.
\end{defn}
By the splitting principle, this is indeed a Pontryagin number of $M$. From its definition it is clear
that $q_m$ vanishes on product manifolds which have no factor of dimension 
at least $4m-4$. (To avoid trivialities, we may assume that $m>1$.)

\begin{prop}
\label{q-number lemma}
Assume that $k$ and $n$ are odd and $\geq 3$. Setting $m=\frac{1}{2}(k+n)$, 
the $q$-number of $X_n^k(c)$ is given by
\begin{align}
\label{q-formula}
    q_m(X_n^k(c))=k\left[\binom{n+k-3}{n}-(k-1)\right]c^n+(n+1)\left[\binom{n+k-3}{n-2}-k\right]c^{n-2} \ .
\end{align}
\end{prop}
\begin{proof}
We compute $q_m(X_n^k(c))$ from the definition using the formula~\eqref{total pontryagin class}:
\begin{align*}
    q_m(X_n^k(c)) &= \Bigl< (n+1)x^2\,\left[ky^{n+k-2}+nx^{n+k-2}+(y+cx)^{n+k-2}\right] \\
    &+ ky^2\,\left[(k-1)y^{n+k-2}+(n+1)x^{n+k-2}+(y+cx)^{n+k-2}\right] \\
    &+ (y+cx)^2\,\left[ky^{n+k-2}+nx^{n+k-2}\right] , [X_n^k(c)]\Bigr> \ .
\end{align*}

The cohomology relations from Lemma~\ref{cohomology} give $x^{n+k-2}=0$ and $(y+cx)\cdot y^k=0$, 
therefore the above sum reduces to
\begin{align*}
    q_m(X_n^k(c))= \Bigl< (n+1)x^2 \left[ky^{n+k-2}+(y+cx)^{n+k-2}\right] , [X_n^k(c)]\Bigr>\\
    + \Bigl< ky^2 \left[(k-1)y^{n+k-2}+(y+cx)^{n+k-2}\right] , [X_n^k(c)] \Bigr> \ .
\end{align*}

We expand the powers of $y+cx$ using the binomial theorem, and, as in the proof of 
Proposition~\ref{Thom number lemma}, we drop terms where the exponent of $x$ is $>n$, 
and we trade powers of $y$ with exponent $>k$ using the relation $y^{k+1}=-cxy^k$ from 
Lemma~\ref{cohomology}. Finally, evaluating on the fundamental class
using $\langle y^kx^n, [X_n^k(c)]\rangle =1$ we find
\begin{align*}
    q_m(X_n^k(c))=(n+1) \left[-k+\sum_{i=0}^{n-2}(-1)^{i+1}\binom{n+k-2}{i}\right]c^{n-2}\\
    +k\left[-(k-1)+\sum_{i=0}^{n}(-1)^{i+1}\binom{n+k-2}{i}\right]c^n \ .
\end{align*} 

 The first sum is simplified using the recursion~\eqref{recursive identity} to get
$$
\sum_{i=0}^{n-2}(-1)^{i+1}\binom{n+k-2}{i}=\binom{n+k-3}{n-2} \ ,
$$
and similarly for the second sum. This completes the proof.
\end{proof}

\begin{rem}\label{linearremark}
When $n=3$, then the shape of formula~\eqref{q-formula} says that $q_m(X_3^k(c))$ is 
an integral linear combination of $c^3$ and $c$. However, an inspection of the linear
term reveals that it actually vanishes in this case, so $q_m(X_3^k(c))$ is simply a 
multiple of $c^3$.
\end{rem}

\section{Calculations in spin bordism}\label{s:bordism}

In this section we work in $\Omega_*$, the rational spin bordism ring, i.e. the usual spin bordism ring tensored with $\Q$.
As is well known, at the rational level there is no difference between oriented bordism and spin bordism, see~\cite[Ch.~XI]{S},
and so $\Omega_*$ is a polynomial algebra over $\Q$ with one ring generator $\alpha_i$ in each dimension of the form $4i$.
Moreover, by a result of Thom~\cite{T}, a $4i$-dimensional closed oriented manifold $M$ may be taken as a representative 
for $\alpha_i$ if and only if its Milnor--Thom number $s_i(M)$ does not vanish; cf.~\cite{HBJ,S}. We will write simply 
$\alpha_i=M\in\Omega_{4i}$, identifying $M$ with its bordism class, without using brackets or other notation to specify that 
we are taking the equivalence class of $M$ in the $\Q$-vector space $\Omega_{4i}$, the degree $4i$ part of $\Omega_*$.

Let $K3$ be the smooth closed oriented $4$-manifold underlying a complex $K3$ surface, and $\HH P^2$ the 
quaternionic projective plane.
\begin{prop}\label{basissequences}
Let $\alpha_1=K3$, $\alpha_2=\HH P^2$. For $i\geq 3$ consider any decomposition $2i=n+k$ with both $n$ and $k$ odd 
and $\geq 3$, and let $\alpha_i(c)=X_n^k(c)$. Then as long as $c$ is even and nonzero, the $\alpha_i$ form a basis
sequence for $\Omega_*$.
\end{prop}
\begin{proof}
The $K3$ surface has trivial first Chern class, and so is spin. Moreover, its signature is not zero, which means $s_1(K3)\neq 0$.
The quaternionic projective plane is $2$-connected, and therefore spin. It is well known that $s_2(\HH P^2)\neq 0$.
Finally, if $k$ and $n$ are odd and $c$ is even, then $ X_n^k(c)$ is spin by Lemma~\ref{spinlemma}, and we have
$s_i(X_n^k(c))\neq 0$ by Proposition~\ref{Thom number lemma} as soon as $c\neq 0$ and $n, k \geq 3$.
\end{proof}

The universal elliptic genus is a surjective homomorphism of graded rings
$$
\varphi\colon\Omega_*\lra\Q[\delta,\epsilon] \ ,
$$
where $\delta$ and $\epsilon$ have degrees $4$ and $8$ respectively; cf.~\cite{HBJ,O}. Ochanine's theorem~\cite{O}
characterises $\ker (\varphi )$ as the ideal generated by all $\C P^{odd}$-bundles with compact Lie groups 
as structure groups. Since the manifolds $X_n^k(c)$ are $\C P^{odd}$-bundles with structure 
group $S^1$, they are in the kernel of the elliptic genus. Moreover, in terms of the above basis sequences we can think 
of $\varphi$ as being the projection 
$$
\Omega_*=\Q[\alpha_1,\alpha_2,\alpha_3(c),\ldots ]\lra\Q[\alpha_1,\alpha_2] 
$$
which kills all the $\alpha_i$ with $i\geq 3$.

\begin{prop}\label{noK3}
Let $p$ and $q$ be odd and $\geq 3$ and $c$ be even and nonzero. Then in $\Omega_*$ the product $X_p^q(c)\times K3$ 
equals a polynomial in $\HH P^2$ and the various $X_n^k(c)$, but not involving $K3$.
\end{prop}
\begin{proof}
The proof is by induction on the dimension $2(p+q)$. The base case is where this dimension is $12$ and $p=q=3$.
Then we look at $X_3^3(c)\times K3$, of dimension $16$. In dimension $16$ we have two choices of indecomposable 
generators, namely $X_3^5(c)$ and $X_5^3(c)$ whose Milnor--Thom numbers according to 
Proposition~\ref{Thom number lemma} are 
$$
s_4(X_3^5(c))=30c^3 \ , \ \ \ \ s_4(X_5^3(c))=18c^5 \ .
$$
This implies that the element $Y(c)=3c^2X_3^5(c)-5X_5^3(c)$ is in the kernel of $s_4$, and is therefore expressible 
as a polynomial in the generators of dimension $\leq 12$. However $Y(c)$ is also in the kernel of the elliptic genus,
and in dimension $16$ this kernel is spanned by $\alpha_3(c)\times K3= X_3^3(c)\times K3$ and $\alpha_4(c)$.
Thus $Y(c)$ is a rational multiple of $X_3^3(c)\times K3$, and we only have to show that it is not the zero multiple.
We can do this conveniently by using the calculations of the $q$-number in Proposition~\ref{q-number lemma}:
\begin{alignat*}{1}
q_4 (Y(c)) &= 3c^2 q(X_3^5(c))-5 q(X_5^3(c))\\
&=3c^2 \cdot 30c^3-5\cdot (-3c^5+42c^3)\\
&=105c^3(c^2-2) \ .
\end{alignat*}
This is nonzero since $c$ is a nonzero integer. Thus we have shown that $X_3^3(c)\times K3$ is a rational
multiple of $Y(c)=3c^2X_3^5(c)-5X_5^3(c)$.

For the inductive step consider some $X_p^q(c)$ of dimension $4m=2(p+q)\geq 16$. We fix a basis sequence
$\alpha_i$ as in Proposition~\ref{basissequences}, taking $X_p^q(c)$ in its dimension: $\alpha_m(c)=X_p^q(c)$
for $m=\frac{1}{2}(p+q)$. Now in dimension $4m+4$ we may consider the two manifolds $X_3^{2m-1}(c)$
and $X_{2m-1}^3(c)$. By Proposition~\ref{Thom number lemma} their Milnor--Thom numbers are
\begin{alignat*}{2}
s_{m+1}(X_3^{2m-1}(c)) &=c^3 \left[\binom{2m+1}{3}-(2m-1)\right] &&= \lambda c^3\ , \\
s_{m+1}(X_{2m-1}^3(c)) &=c^{2m-1} \left[\binom{2m+1}{2m-1}-3\right]  &&= \mu c^{2m-1}\ .
\end{alignat*}
This implies that the element $Z(c)=\mu c^{2m-4}X_3^{2m-1}(c)-\lambda X_{2m-1}^3(c)$ is in the kernel of 
$s_{m+1}$, and is therefore expressible as a polynomial in the generators $\alpha_i$ with  $i\leq m$.

We want to check that $q_{m+1}(Z(c))\neq 0$. Using the formula from Proposition~\ref{q-number lemma}
we see that $q_{m+1}(Z(c))$ is an integral linear combination of $c^{2m-1}$ and $c^{2m-3}$. However,
in light of Remark~\ref{linearremark}, $c^{2m-4}X_3^{2m-1}(c)$ does not contribute to the term of 
degree $2m-3$, and so we can simply read off the coefficient of $c^{2m-3}$ in $q_{m+1}(Z(c))$ to be
$$
-\lambda\cdot 2m \cdot \left[\binom{2m-1}{2m-3}-3\right] 
= -\left[\binom{2m+1}{3}-(2m-1)\right] \cdot 2m \cdot \left[\binom{2m-1}{2}-3\right] \ .
$$
As this is clearly nonzero, we do indeed have $q_{m+1}(Z(c))\neq 0$ for nonzero $c$.

Now expressing $Z(c)$ as a polynomial in the generators $\alpha_i$ with  $i\leq m$, the non-vanishing 
of the $q$-number tells us that the monomial $\alpha_1\cdot\alpha_m(c)$ must appear with a nonzero
coefficient, since $q_{m+1}$ vanishes on all monomials consisting only of $\alpha_i$ with $i<m$. 
We can then solve the resulting equation for $\alpha_1\cdot\alpha_m(c)$ and this expresses 
$K3\times X_p^q(c)$ as a rational linear combination of $X_3^{2m-1}(c)$, $X_{2m-1}^3(c)$ and 
monomials in the $\alpha_i$ with  $i< m$.
As $K3\times X_p^q(c)$ is in the kernel of the elliptic genus, each monomial appearing in the 
linear combination must contain an $\alpha_i$ with $3\leq i<m$ represented by some $X_n^k(c)$. 
By the inductive hypothesis, whenever such a monomial also contains $\alpha_1=K3$ we can 
replace $K3\times X_n^k(c)$ by an expression not involving $K3$. This completes the inductive step.
\end{proof}

After these preparations we can now prove the theorem.

\begin{proof}[Proof of Theorem~\ref{t:main}]
We fix a basis sequence $\alpha_1=K3$, $\alpha_2=\HH P^2$, $\alpha_i(c)=X_n^k(c)$ for $i\geq 3$ 
for $\Omega_*$ as in Proposition~\ref{basissequences}. As discussed above, the elliptic genus $\varphi$
is the projection to the quotient by the ideal generated by the $\alpha_i$ with $i \geq 3$.

Consider spin manifolds of dimension $4m$, and a non-trivial linear combination $f$ of their Pontryagin numbers.
If the linear map 
$$
f\colon \Omega_{4m}\lra \Q
$$ 
does not factor through $\varphi\vert\Omega_{4m}$, then it does not vanish identically on $\ker (\varphi)$. 
This means that $m\geq 3$ and that there is an element in $\ker (\varphi)$ on which $f$ does not vanish. However,
a vector space basis for $\ker (\varphi)$  is given by the monomials in the $\alpha_i$ containing at least
one index $i\geq 3$. Thus there is a Cartesian product of $K3$, $\HH P^2$
and the $X_n^k(c)$ on which $f$ does not vanish, and this product definitely contains at least 
one factor of the form $X_n^k(c)$. If this product does not contain $K3$, then it is a product 
of nonnegatively curved manifolds, and therefore nonnegatively curved. The value of $f$ on this 
Cartesian product can be thought of as a linear combination of Pontryagin numbers of one of the 
factors of the form $X_n^k(c)$, with coefficients that depend on the Pontryagin numbers of the other factors. 
Moreover, the non-vanishing of $f$ means that a non-zero Pontryagin number of $X_n^k(c)$ does appear. 
Since the Pontryagin numbers of $X_n^k(c)$ are odd polynomials in $c$ by Proposition~\ref{prop:odd}, 
they are certainly non-constant. Varying $c$, we see that $f$ is unbounded on this family of nonnegatively curved 
manifolds.

Finally, if the monomial in the generators $\alpha_i$ on which $f$ does not vanish contains a $K3$ factor,
 then, because it also contains at least one $X_n^k(c)$, we can replace $K3\times X_n^k(c)$ by a 
 linear combination of terms not involving $K3$ using Proposition~\ref{noK3}. 
 This replaces the monomial in the generators by a linear combination of terms, but since $f$ does not
 vanish on the monomial, it does not vanish on at least one of the summands of the linear combination.
 Repeating this procedure until there are 
 no $K3$ factors left we find a product of factors all of which are $\HH P^2$ or of the form $X_n^k(c)$,
 and on which $f$ does not vanish. Since there is at least one $X_n^k(c)$-factor, as above 
 we conclude that $f$ is unbounded on these nonnegatively curved manifolds as we vary $c$.
 This completes the proof.
\end{proof}

\bigskip

\bibliographystyle{amsplain}

\bigskip

\end{document}